\documentclass[a4paper,12pt]{amsart}
\usepackage[margin=3cm]{geometry}
\usepackage{amssymb,amsmath}
\newtheorem{theorem}{Theorem}[section]

\newtheorem{proposition} [theorem]{Proposition}
\newtheorem{lemma}[theorem]{Lemma}
\newtheorem{definition}[theorem]{Definition}

\newcommand{\dquer}{\overline\partial}
\newcommand{\dquers}{\overline\partial ^*_\varphi}
\newcommand{\boxphi}{\square_\varphi}
\newcommand{\levim}{\frac{\partial^2\varphi}{\partial z_j\partial\overline z_k}}

\numberwithin{equation}{section}
\begin{document}
\title{On the weighted $\dquer $-Neumann problem  on unbounded domains.}

\author{Klaus Gansberger}

\thanks{Supported by the FWF-grant  P19147.}

\address{K. Gansberger, Institut f\"ur Mathematik, Universit\"at Wien,
Nordbergstrasse 15, A-1090 Wien, Austria.}
\email{klaus.gansberger@univie.ac.at}

\maketitle

\begin{abstract} ~\\
Let $\Omega$ be an unbounded, pseudoconvex domain in $\Bbb C^n$ and let $\varphi$ be a $\mathcal C^2$-weight function plurisubharmonic on $\Omega$.  We show both necessary and sufficient conditions for existence and compactness of a weighted $\dquer$-Neumann operator $N_\varphi$ on the space $L^2_{(0,1)}(\Omega,e^{-\varphi})$ in terms of the eigenvalues of the complex Hessian $(\partial ^2\varphi/\partial z_j\partial\overline z _k)_{j,k}$ of the weight. We also give some applications to the unweighted $\dquer$-Neumann problem on unbounded domains.
\end{abstract}

\section{Introduction.}~\\

The subject of this paper is the weighted $\dquer$-Neumann problem on pseudoconvex, unbounded domains. The weighted $\dquer$-Neumann operator is the inverse of the weighted complex Laplacian, which acts on (p,q)-forms that satisfy certain boundary conditions, see Section \ref{pre} for the precise definitions. The weighted $\dquer$-equation is one of the fundamental tools in complex analysis, see e.g. \cite{h}. It also arises when studying the unweighted problem: For instance in the case of complete pseudoconvex Hartogs domains, the $\dquer$-Neumann problem can be reduced to a corresponding weighted problem on the base domain  \cite{be}, \cite{li}. A third motivation comes from the study of three-dimensional, pseudoconvex, compact CR-manifolds, see \cite{ch}. 
\vskip 0.5 cm 

The unweighted $\dquer$-Neumann problem on bounded domains has been intensively studied and is of interest in complex analysis for various reasons. For background on the $\dquer$-Neumann problem, we refer the reader to \cite{bs}, \cite{fk}  and \cite{ChSh}.\\
One reason for the interest in this problem is that existence of a bounded $\dquer$-Neumann operator implies solvability of the inhomogeneous $\dquer$-equation with control of the norm of the solution (a priori only in the $L^2$-sense).  The question of compactness of $N$ is of interest for its own right, see for instance \cite{FS2} for a discussion. To mention one of the most important reasons, compactness of $N$ implies global regularity in the sense of preservation of Sobolev spaces, see \cite{KN}. This in turn has consequences for the extension behavior of biholomorphisms.\\
More recently, compactness is being studied not only as a property stronger than global regularity, but also as one for which a characterization in terms of the boundary should be possible, whereas global regularity seems to be too subtle and unstable for this. Generally, compactness is  believed to be more tractable than global regularity. 
\vskip 0.5 cm
 
In \cite{Ca}, Catlin introduced his notion of Property (P), giving a sufficient condition for compactness of $N$ and thus for global regularity, which can be verified on a large class of domains. A bounded, smooth domain $\Omega$ is said to satisfy Property (P), if for each $M\in\mathbb N$ there is  a function $\lambda_M\in\mathcal C^\infty(\overline\Omega)$, such that $0\le\lambda_M\le1$ and for all $p\in\partial \Omega$ and all $t\in\mathbb C^n$
$$\sum_{j,k=1}^n\frac{\partial^2\lambda_M}{\partial z_j\partial\overline z_k}(p)t_j\overline t_k\ge M\|t\|^2.$$
McNeal gave a generalization -- Property (\~ P) -- still implying compactness, see \cite{mc}. He replaced the uniform boundedness of the family by self-boundedness of the complex gradient, \cite{mc} Definition 1. One can easily check that Property (P) always implies Property (\~ P). There are some cases known, in which Property (P) turns out to be also necessary for compactness of $N$ (see e.g. \cite{fs}), but in general it is not understood how much room there is between compactness and Property (P) or between Property (\~ P) and Property (P).
\vskip 0.5 cm 

Few is known for the case of unbounded domains. Recent contributions to the $\dquer$-Neumann problem in weighted $L^2$-spaces on $\mathbb C$ are \cite{hahe} and \cite{ortega}.  Weighted spaces on $\mathbb C^n$ were considered in \cite{hahe} and \cite{gh}. In the present paper we develop methods used in \cite{gh} further which allows us to also treat unbounded pseudoconvex domains with boundary.\\
The main result on existence is the following.
\begin{theorem}
\label{existence}
Let $\Omega$ be a smooth, pseudoconvex, unbounded domain and denote by $\lambda_\varphi(z)$ the lowest eigenvalue of the complex Hessian $(\partial ^2\varphi/\partial z_j\partial\overline z _k)_{j,k}$ of the weight function. Suppose that 
\begin{equation}
\liminf_{z\in\Omega, |z|\to\infty}\lambda_\varphi\ge\varepsilon
\end{equation}
for some $\varepsilon>0$. Then there exists a bounded $\dquer$-Neumann operator on $L^2(\Omega, \varphi)$.
\end{theorem}
To formulate the sufficient condition for compactness, we need a notion of Property (P) for unbounded domains. We shall use the following local version.

\begin{definition}
\label{prop p}
An unbounded, smooth domain satisfies Property (P), if  the following holds: for any $p\in\partial \Omega$ there is a neighborhood $U_p$ such that for each $M\in\mathbb N$ there is  a function $\varphi_{p,M}\in PSH(U_p)\cap\mathcal C^\infty(\overline U_p)$, with $0\le\varphi_{p,M}\le 1$ and $\lambda_{\varphi_{p,M}}\ge M $ on $U_p\cap \partial\Omega$, where $\lambda_{\varphi_{p,M}}\ge M$ denotes the lowest eigenvalue of the complex Hessian $(\partial ^2\varphi_{p,M}/\partial z_j\partial\overline z _k)_{j,k}$.
\end{definition}

If the domain is bounded, this Definition coincides with the original one of Catlin, as it was stated above. \\

\begin{theorem}
\label{cp 3}
Let $\Omega$ be unbounded, smooth and pseudoconvex and suppose that its boundary satisfies Property (P). Let $\lambda_\varphi(z)$ denote the lowest eigenvalue of the complex Hessian $(\partial ^2\varphi/\partial z_j\partial\overline z _k)_{j,k}$ of the weight function. Suppose furthermore that 
$$\lim_{|z|\to\infty, z\in\Omega}\lambda_\varphi(z)=\infty.$$ 
Then the weighted $\dquer$-Neumann operator $N_\varphi$ exists and is a compact operator from $L^2_{(0,1)}(\Omega,\varphi)$ into itself. 
\end{theorem}

To also give a necessary condition, we prepare the following Definition taken from \cite{adams}.

\begin{definition}
We call a domain $\Omega$ quasibounded if and only if $$\lim_{z\in\Omega, \ |z|\to\infty} dist(z,\partial\Omega)=0.$$
Equivalently, $\Omega$ is quasibounded if and only if there is no $r>0$ such that $\Omega$ contains a sequence of congruent pairwise disjoint balls with radius $r$.
\end{definition}

\textbf{Remark. } Although a general quasibounded domain can be much more complicated, one can typically think of such a domain to look like $\Omega=\{(z,w)\in\Bbb C^2\ :\ |zw|<1\}$. For further details on the notion of quasiboundedness, see \cite{adams}, Chapter 6.

\begin{theorem}
\label{ness}
Suppose that $\Omega$ is an unbounded but not quasibounded domain and suppose that $N_\varphi$ is a compact operator on $L^2_{(0,1)}(\Omega, \varphi)$. Then for any sequence $\mathbb B (z_l,r)$  of disjoint balls with fixed radius $r$ contained in $\Omega$ it holds
$$\lim_{l\to\infty}\int_{\mathbb B (z_l,r)}\triangle\varphi\ d\lambda=\infty.$$
\end{theorem}

\textbf{Remark. } For plurisubharmonic functions, $\triangle\varphi$ is comparable to the largest eigenvalue $\lambda_n$ of the complex Hessian $(\partial ^2\varphi/\partial z_j\partial\overline z _k)_{j,k}$. Thus, one can think of $\int_{\mathbb B (z_l,r)}\triangle\varphi\ d\lambda$ to be a regularized version of $\lambda_n$. Theorem \ref{ness} states that compactness of $N_\varphi$ implies that the mean value of $\lambda_n$ has to tend to infinty at infinity, which should be compared with the condition from Theorem \ref{cp 3}.
\vskip 0.5 cm

\textbf{Remark. } For the case $\Omega=\mathbb C$, it was shown in \cite{hahe} that  
$$\lim_{l\to\infty}\int_{\mathbb B (z_l,r)}\left(\triangle\varphi\right)^2d\lambda=\infty$$
for any sequence $(\mathbb B (z_l,r))_l$ of  disjoint balls with $|z_l|\to\infty$ is necessary and sufficient for compactness if one assumes $\triangle\varphi\in B_2$, a reverse H\"older class. In fact, this condition is necessary for compactness for $\Omega=\mathbb C^n$, $n\ge1$, and arbitrary plurisubharmonic weight function $\varphi$, as was shown in \cite{ga}. Both \cite{hahe} and \cite{ga} apply spectral analytic Theorems to prove that result -- in contrast to the more direct and purely complex analytic we give here, yielding a sharper result.\\
Marzo and Ortega-Cerd\'{a} showed in \cite{ortega} under the condition that $\mu=\triangle\varphi\ d\lambda$ defines a doubling measure that
$$\lim_{l\to\infty}\int_{\mathbb B (z_l,r)}\triangle\varphi\ d\lambda=\infty$$
is equivalent to compactness of the canonical solution operator $\dquers N_\varphi$ to $\dquer$ in $L^2(\Bbb C, \varphi)$. This was done by carefully estimating the Bergman kernel.
\vskip 1cm

\section{Preliminaries.}~\\
\label{pre}

Let $\Omega$ an unbounded pseudoconvex domain in $\Bbb C^n$ with smooth boundary, i.e., there is a smooth function $r:\Bbb C^n\to\Bbb R$ such that $\Omega=\{ z\in\Bbb C^n\ |\ r(z)<0\}$ with $|\nabla r|\neq0$ on the set $\{r=0\}$ and 
$$\sum_{j,k=1}^n\frac{\partial^2r}{\partial z_j\partial \overline z_k}(p)t_j\overline t_k\ge 0 $$
for all $p\in\partial\Omega$ and all $t\in T_p^{0,1}\partial\Omega$. Let furthermore $\varphi : \Omega\longrightarrow \mathbb R^+ $ be a plurisubharmonic weight function of class $\mathcal C^2$ and define the space
$$L^2(\Omega , \varphi )=\{ f:\Omega \longrightarrow \mathbb C \ | \ \int_{\Omega}
|f|^2\, e^{-\varphi}\,d\lambda < \infty \},$$
where $\lambda$ denotes the Lebesgue measure. Similarly define the space $L^2_{(0,1)}(\Omega, \varphi )$ of $(0,1)$-forms with coefficients in $L^2(\Omega , \varphi )$ and the space $L^2_{(0,2)}(\Omega, \varphi )$ of $(0,2)$-forms with coefficients in
$L^2(\Omega , \varphi ).$
Let 
$$\langle f,g\rangle_\varphi=\int_{\Omega}f \,\overline{g} e^{-\varphi}\,d\lambda$$
denote the inner product and 
$$\| f\|^2_\varphi =\int_{\Omega}|f|^2e^{-\varphi}\,d\lambda $$
the norm in $L^2(\Omega , \varphi ).$ Defining the $\dquer$-operator, we set on $\mathcal C_0^\infty(\Omega)$, i.e. the space of smooth functions with compact support in $\Omega$, 
$$\dquer f =\sum_{j=1}^n\frac{\partial f}{\partial \overline z_j}d\overline z_j.$$
Taking the maximal closure of this operator and still denoting it by $\dquer,$ we turn $\dquer$ into a closed, densely defined operator on $L^2(\Omega, \varphi )$. Moreover, it can be extended to $(0,q)$-forms in the natural way by setting
$$\dquer f= \sum_{j,K}\frac{\partial f_K}{\partial \overline z_j}d\overline z_j\wedge d\overline z_K$$
for $f=\sum_{|K|=q}f_kd\overline z_K$. As a closed, densely defined operator, $\dquer$ possesses a Hilbert space adjoint which we denote by $\dquers$. For $f=\sum_{j=1}^nf_jd\overline z_j\in dom (\dquer_\varphi^*)$ one has
$$\dquer_\varphi^*f=-\sum_{j=1}^n \left ( \frac{\partial}{\partial z_j}-
\frac{\partial \varphi}{\partial z_j}\right )f_{j}.$$
The complex Laplacian on $(0,1)$-forms is defined to be
$$\boxphi = \dquer  \,\dquers + \dquers \dquer.$$
This is a closed, selfadjoint and positive operator, which means that 
$$\langle \boxphi f,f\rangle_\varphi \ge 0 \ , \   {\text{for}} \  f\in dom (\boxphi ).$$
The associated Dirichlet form is 
$$Q_\varphi (f,g)= \langle \dquer f,\dquer g\rangle_\varphi + \langle \dquers f ,\dquers g\rangle_\varphi, $$
with form domain $dom (\dquer ) \cap dom (\dquers ).$ The weighted $\dquer $-Neumann operator on the level of $(0,1)$-forms, which we denote by 
$N_\varphi $, is -- if it exists -- the bounded inverse of $\boxphi .$ Note that we see by the same argument as in  \cite{gh}, Lemma 2.3, that existence and compactness of $N_\varphi$ is invariant under equivalent weights, where we call two weight functions equivalent if the weighted $L^2$-norms induced are equivalent. Thus without loss of generality, we restrict ourselves from now on to smooth weight functions.
\vskip 1 cm

\section{The weighted problem.}~\\

To begin with, let us give the simple characterization of the domain of $\dquers$  in the weighted space $L^2_{(0,1)}(\Omega, \varphi)$.

\begin{proposition}
\label{dom dquers}
Let $f=\sum f_jd\overline z_j\in L^2_{(0,1)}(\Omega,\varphi)$ and let $r$ be a defining function of $\Omega$ with $|\nabla r|^2=1$ on $\partial\Omega$. Then $f\in dom(\dquers)$ if and only if $\sum _{j=1}^nf_j\frac{\partial r}{\partial z_j}=0$ on $\partial\Omega$ as well as
\begin{equation*}
 \sum_{j=1}^n\left ( \frac{\partial f_j}{\partial z_j}- \frac{\partial \varphi}{\partial z_j}\, f_j \right )
 \in L^2(\Omega, \varphi ).
\end{equation*}
\end{proposition}

{\em Proof.}
Let a function $f$ fulfilling the conditions be given and let $(\chi_R)_{R\in\mathbb N}$ be a family of smooth cutoff functions identically one on $\mathbb B_R$, the ball with radius $R$, and supported in $\mathbb B_{R+1}$. Suppose additionally that all first order derivatives of the functions in this family are uniformly bounded by a constant $M$. Then for all $g\in dom(\dquer)$ we have via integration by parts
\begin{align*}
\langle \chi _Rf,\dquer g\rangle_\varphi=&-\int\limits_{\Omega}\sum_{j=1}^n \frac{\partial}{\partial  z_j}\left(\chi_R f_je^{-\varphi}\right)\overline g\ d\lambda+\int\limits_{\partial\Omega}\chi_R\sum_{j=1}^nf_j\frac{\partial r}{\partial z_j}e^{-\varphi}\overline g\ d\sigma\\
=&-\int\limits_{\Omega}\sum_{j=1}^n \frac{\partial}{\partial  z_j}\left(\chi_R f_je^{-\varphi}\right)\overline g\ d\lambda.
\end{align*}
Now doing the limit $R\to\infty$, it is easily seen that
$$|\langle f,\dquer g\rangle_\varphi|\le \|g\|_\varphi\left\Vert e^\varphi\sum_{j=1}^n\frac{\partial}{\partial z_j}\left( f_je^{-\varphi}\right)\right\Vert_\varphi+M\Vert g\Vert_\varphi\Vert f\Vert_\varphi.$$
So by assumption $|\langle f,\dquer g\rangle_\varphi|\leq C\Vert g\Vert_\varphi$ and thus $f\in dom(\dquers)$. Conversely, for $f\in dom(\dquer^*_\varphi)$ and any $g\in \mathcal C_0^\infty(\Omega)$,
 \begin{align*}
 \langle \dquer^*_\varphi f,g \rangle_\varphi=&\langle f,\dquer g\rangle_\varphi
 =\sum_{j=1}^n\ \int\limits_\Omega f_j\frac{\partial\overline g}{\partial z_j}e^{-\varphi}d\lambda.\\
  \end{align*}
Since $\mathcal C_0^\infty(\Omega)$ is dense in $L^2(\Omega, \varphi )$, we get after integrating by parts that
$$\dquer^*_\varphi f=-e^{\varphi}\sum_{j=1}^n\frac{\partial}{\partial z_j}\left( fe^{-\varphi}\right),$$
hence in particular $e^{\varphi}\sum_{j=1}^n\frac{\partial}{\partial z_j}\left( f_je^{-\varphi}\right)\in L^2_\varphi(\Omega)$. Doing the same calculation for general $g\in dom(\dquer)$, integration by parts again yields
$$ \langle \dquer^*_\varphi f,g \rangle_\varphi=\langle f,\dquer g\rangle_\varphi=\langle -e^{\varphi}\sum_{j=1}^n\frac{\partial}{\partial z_j}\left( fe^{-\varphi}\right),g\rangle_\varphi+\int\limits_{\partial\Omega}g\sum_{j=1}^nf_j\frac{\partial r}{\partial z_j}e^{-\varphi}d\sigma.$$
Thus by comparing the two expressions for $\dquers f$, we see that the boundary integral has to vanish for all $g$, which is the case if and only if $\sum _{j=1}^nf_j\frac{\partial r}{\partial z_j}=0$ on $\partial\Omega$.
\begin{flushright}
$\square$
\end{flushright}
 
The following Lemma generalizes a well-known density Lemma to unbounded domains and is the first important technical step in our considerations.

\begin{lemma}
\label{density}
Let $r$ be a defining function of $\Omega$ such that $|\nabla r|^2=1$ on $\partial\Omega$ and suppose that $\partial \Omega$ is of class $\mathcal C^{k+1}$. Then for any $f\in dom(\dquer)\cap dom(\dquers)$ there is a sequence $(f^{(l)})_l\subset \mathcal C_{(0,1)}^k(\overline\Omega)$ such that  $f^{(l)}\to f$ in the graph norm $f\mapsto (\Vert f\Vert^2 _\varphi+\Vert \dquer f\Vert^2 _\varphi+\Vert \dquers f\Vert ^2_\varphi)^\frac{1}{2}$ and $f^{(l)}$ vanishes on $\Omega\setminus\mathbb B_{l+1}$ as well as $\sum _{j=1}^nf^{(l)}_j\frac{\partial r}{\partial z_j}$ on $\partial\Omega$.
\end{lemma}

{\em Proof.}
Keeping the notation from Lemma \ref{dom dquers}, we easily see by a direct computation that $\chi _l f\to f$ in the graph norm as $l\to\infty$. Now using Lemma 4.3.2 in \cite{ChSh} for each fixed $l$, the function $\chi_l f$ can be approximated by a sequence of functions with the claimed smoothness properties, fulfilling the boundary condition and support in $\mathbb B_{l+1}$. Thus the Lemma follows  by choosing an appropriate diagonal sequence.
\begin{flushright}
$\square$
\end{flushright}

\begin{proposition}
\label{kohn-morrey}
\textbf{(Kohn -- Morrey formula)}
Let $\Omega$ be of class $\mathcal C^2$  and let $r$ be a defining function of $\Omega$ such that $|\nabla r|^2=1$ on $\partial\Omega$. Then for any $f=\sum_{j=1}^n f_jd\overline z_j\in dom(\dquer)\cap dom(\dquers)$
\begin{align*}
\sum_{j,k=1}^n\ \int \limits_{\Omega}\levim f_j\overline f_ke^{-\varphi}d\lambda&+\sum_{j,k=1}^n\  \int \limits_{\Omega}\left| \frac{\partial f_j}{\partial \overline z_k}\right| ^2e^{-\varphi}d\lambda+\sum_{j,k=1}^n \ \int\limits_{\partial\Omega}\frac{\partial^2 r}{\partial z_j\partial\overline z_k}f_j\overline f_ke^{-\varphi}d\sigma\\
&=\Vert\dquer f\Vert_\varphi^2+\Vert \dquer^*_\varphi f\Vert^2_\varphi ,
\end{align*}
where $\sigma$ denotes the surface measure on $\partial \Omega$.
\end{proposition}
{\em Proof.}
Using Lemma \ref{density}, the Proposition follows from the Kohn -- Morrey formula on bounded domains by the Dominated Convergence Theorem.\\
See for instance  \cite{ChSh}, Proposition 4.3.1, for a proof in the bounded case.
\begin{flushright}
$\square$
\end{flushright}

From this identity we can immediately conclude Theorem \ref{existence}.\\

{\em Proof of Theorem \ref{existence}.}
Since existence of $N_\varphi$ is invariant under equivalent weights, we can after possibly shrinking $\varepsilon$ without loss of generality assume that $\lambda_\varphi(z)\ge\varepsilon$ for all $z\in\Omega$. Since $\Omega$ is pseudoconvex, Proposition \ref{kohn-morrey} yields $\varepsilon\| f\|_\varphi\le \Vert\dquer f\Vert_\varphi^2+\Vert \dquer^*_\varphi f\Vert^2_\varphi $ for all $f\in dom(\dquer)\cap dom(\dquers)$, so $N_\varphi$ is bounded. 
\begin{flushright}
$\square$
\end{flushright}
\vskip 1 cm

\section{Weighted Sobolev spaces}~\\

Similar to the case of bounded domains, our strategy to find a sufficient condition for compactness of the weighted $\dquer$-Neumann operator $N_\varphi$ is to show a so-called compactness estimate (see Proposition \ref{compact}). To this end, we need a norm on $L^2(\Omega,\varphi)$ that is strictly weaker than the weighted $L^2$-norm. On bounded domains, one naturally has the Sobolev norm $\| . \|_{-1}$, which is strictly weaker than the $L^2$-norm by the Rellich -- Kondrachov Theorem. On unbounded domains, it is in general not true that $H^1(\Omega)$ embeds compactly into $L^2(\Omega)$. Thus we need an appropriate notion of a weighted Sobolev space and a compact injection into $L^2(\Omega,\varphi)$. Similar Definitions in fact already appeared before in \cite{gh}.

\begin{definition}
\label{sob}
Denote the coordinates in $\Bbb C^n$ by $(z_1,\dots ,z_n)=(x_1,y_1,\dots ,x_n,y_n)$. For $k\in \mathbb{N}$ let 
$$H^k(\Omega, \varphi):= \{ f\in L^2(\Omega, \varphi) \ | \ D^\alpha f \in L^2(\Omega, \varphi) \ {\text{for any}} \ |\alpha | \le k \},$$
where $D^\alpha =\frac{\partial^{| \alpha|}}{\partial^{\alpha_1}x_1\dotsm \partial^{ \alpha_{2n}}y_n}$, with the norm
$$\|f\|_{k,\varphi}^2= \sum_{|\alpha |\le k}\|D^\alpha f\|_\varphi ^2.$$
Let moreover $H^k_0(\Omega, \varphi)$ be the closure of $\mathcal C_0^\infty(\Omega) $ under the norm defined above.
\end{definition}

\begin{definition}
\label{dualsob}
For $j=1,\dots ,n$ let
$$X_j=\frac{\partial }{\partial x_j} - \frac{\partial \varphi}{\partial x_j} \ {\text {and}} \ 
Y_j=\frac{\partial }{\partial y_j} - \frac{\partial \varphi}{\partial y_j},$$
and define
$$H^k(\Omega, \varphi ,\nabla \varphi)=\{ f\in L^2(\Omega, \varphi) \ | \ T^\alpha f\ \in  L^2(\Omega, \varphi) ,\ {\text{for any}} \ |\alpha | \le k \},$$
where $T^\alpha=X_1^{\alpha_1}Y_1^{\alpha_2}\dotsm X_n^{\alpha_{2n-1}}Y_n^{\alpha_{2n}}$, with the norm
$$\|f\|^2_{1, \varphi ,\nabla \varphi}= \sum_{|\alpha |\le k}\|D^\alpha f\|_\varphi ^2 .$$
Similarly, define $H_0^k(\Omega, \varphi ,\nabla \varphi)$ to be the closure of  $\mathcal C_0^\infty(\Omega) $ under the norm above.
\end{definition}

Note that $X_j$ is the formal adjoint of $-\frac{\partial }{\partial x_j}$ with respect to the weighted inner product. It holds that $X_j^*=-\frac{\partial }{\partial x_j}=-D_j$. The two norms defined above are related in the following way.

\begin{lemma}
\label{sob}
Let $\varphi$ be a plurisubharmonic weight function. Then for any $f\in H_0^1(\Omega, \varphi ,\nabla \varphi)$:
\begin{enumerate}
\item $\|f\|^2_{1, \varphi }\le\|f\|^2_{1, \varphi ,\nabla \varphi}$
\item $\| (D_j\varphi) f\|^2_\varphi\le2 \|f\|^2_{1, \varphi ,\nabla \varphi}$, where $D_j=\frac{\partial }{\partial x_j}$.
\end{enumerate}
\end{lemma}

{\em Proof.}
For any $f\in \mathcal C_0^\infty(\Omega)$ we have $(X_j+X^*_j)f=-\frac{\partial \varphi}{\partial x_j}\, f $ and $
[X_j,X^*_j]f= -\frac{\partial^2\varphi}{\partial x_j^2}\, f$. Thus
\begin{align*}
\|f\|_{1,\varphi}^2= &\| f\|^2_\varphi+\sum_{j=1}^n(\|X^*_jf\|_\varphi ^2+\|Y^*_jf\|_\varphi ^2)\\
=&\| f\|^2_\varphi+\sum_{j=1}^n(\|X_jf\|_\varphi ^2+\|Y_jf\|_\varphi ^2) - \langle\triangle\varphi f,f\rangle_\varphi\\
\le&\| f\|^2_\varphi+\sum_{j=1}^n(\|X_jf\|_\varphi ^2+\|Y_jf\|_\varphi ^2)\\
=&\|f\|^2_{1, \varphi ,\nabla \varphi}.
\end{align*}
By density of $\mathcal C_0^\infty(\Omega)$, this holds for all $f\in H_0^1(\Omega, \varphi ,\nabla \varphi)$.
Now if $D_j=\frac{\partial}{\partial x_j}$, then 
$$\|\varphi_{x_j}f\|^2_\varphi=\|(X_j+X_j^*)f\|^2_\varphi\le\|X_jf\|_\varphi ^2+\|X^*_jf\|_\varphi ^2\le2\|f\|^2_{1, \varphi ,\nabla \varphi}$$
\begin{flushright}
$\square$
\end{flushright}

\textbf{Remark. }On bounded domains, these two Definitions coincide with the classical Definition of a Sobolev space, if one assumes the weight function to be smooth on $\overline \Omega$. Even on unbounded domains, they are equivalent to the usual one if the weight and its first order derivatives are bounded, in particular if the weight is zero. In this sense, the Definitions \ref{dualsob} and \ref{sob} extend the common notion of a Sobolev space.\\
Moreover, $H_0^k(\Omega, \varphi ,\nabla \varphi)\hookrightarrow H^k_0(\Omega, \varphi)$ continuously by Lemma \ref{sob} and thus also $H^{-k}_0(\Omega, \varphi)\hookrightarrow H_0^{-k}(\Omega, \varphi ,\nabla \varphi)$, where we use the convention to denote the dual space of $H_0^k(\Omega, \varphi ,\nabla \varphi)$ by $H_0^{-k}(\Omega, \varphi ,\nabla \varphi)$.
\vskip 0.5 cm

\begin{lemma}
\label{vergl}
Let $\Omega$ be a domain in $\Bbb C^n$ and let $\varphi$ be a $\mathcal C^2$-function. Suppose that $\chi$ is a smooth function with compact support in $\Omega$. Then 
$$\int_\Omega \triangle\varphi\ \chi^2d\lambda=-4\|\nabla\chi\|^2_{L^2}+\int_\Omega\chi^2\ |\nabla\varphi|^2d\lambda .$$
\end{lemma}
{\em Proof.}
In the proof of Lemma \ref{sob} we had the identity
\begin{equation}
\label{sob1}
\|f\|^2_{1, \varphi ,\nabla \varphi}=\|f\|_{1,\varphi}^2 +\langle\triangle\varphi f,f\rangle_\varphi.
\end{equation}
Thus, taking $f=\chi e^{\varphi /2}$, we obtain
$$\|f\|^2_{1, \varphi ,\nabla \varphi}=\sum_{j=1}^n\int_\Omega\left(|\chi_{x_j}-\frac{1}{2}\chi\varphi_{x_j}|^2+|\chi_{y_j}-\frac{1}{2}\chi\varphi_{y_j}|^2\right)d\lambda$$
and 
$$\|f\|_{1,\varphi}^2=\sum_{j=1}^n\int_\Omega\left(|\chi_{x_j}+\frac{1}{2}\chi\varphi_{x_j}|^2+|\chi_{y_j}+\frac{1}{2}\chi\varphi_{y_j}|^2\right)d\lambda.$$
Since the integrand is real-valued, plugging this into \eqref{sob1} it follows by elementary algebra that the difference $\|f\|^2_{1, \varphi ,\nabla \varphi}-\|f\|_{1,\varphi}^2$ equals
\begin{align*}
&\langle\triangle\varphi f,f\rangle_\varphi=\\
&-4\sum_{j=1}^n\int_\Omega\left((\chi_{x_j}+\frac{1}{2}\chi\varphi_{x_j})(\chi_{x_j}-\frac{1}{2}\chi\varphi_{x_j})+(\chi_{y_j}+\frac{1}{2}\chi\varphi_{y_j})(\chi_{y_j}-\frac{1}{2}\chi\varphi_{y_j})\right)d\lambda,
\end{align*}
which implies the Lemma.
\begin{flushright}
$\square$
\end{flushright}

\begin{proposition}
\label{rellich}
Suppose that the weight function satisfies
\begin{align*}
&\lim_{z\in\Omega,|z|\to \infty}(\theta |\nabla \varphi (z)|^2+\triangle \varphi (z))= +\infty \ \ \text{as well as}\\
&\lim_{z\in\Omega,z\to \partial\Omega}(\theta |\nabla \varphi (z)|^2+\triangle \varphi (z))= +\infty
\end{align*}
for some $\theta \in (0,1)$. Then the embedding of $H_0^1(\Omega, \varphi ,\nabla \varphi)$ into $ L^2(\Omega, \varphi) $ is compact.
\end{proposition}

{\em Proof.}
As noted above, for the vector fields $X_j$  and their formal adjoints $X_j^*=-\frac{\partial}{\partial x_j}$ the following relations hold on $  \mathcal{C}^\infty_0(\Omega)$:
$$(X_j+X^*_j)f=-\frac{\partial \varphi}{\partial x_j}\, f \ {\text{and}} \ 
[X_j,X^*_j]f= -\frac{\partial^2\varphi}{\partial x_j^2}\, f,$$
as well as
$$\langle [X_j,X^*_j]f,f\rangle_\varphi=\|X^*_jf\|^2_\varphi - \|X_jf\|^2_\varphi ,$$
$$\|(X_j+X^*_j)f\|^2_\varphi \le (1+1/\epsilon)\|X_jf\|^2_\varphi + (1+\epsilon)\|X^*_jf\|^2_\varphi ,$$
for each $\varepsilon>0$, and similarly for the vector fields $Y_j.$
It follows that
$$\langle |\nabla \varphi (z)|^2+(1+\epsilon)\triangle \varphi (z)f,f\rangle_\varphi \le 
( 2+\epsilon +1/\epsilon)\sum_{j=1}^n( \|X_jf\|^2_\varphi
+ \|Y_jf\|^2_\varphi) ,$$
and since $\mathcal{C}^\infty_0(\Omega)$ is dense in $H_0^1(\Omega, \varphi ,\nabla \varphi)$ by Definition, this inequality is valid for all $f\in H_0^1(\Omega, \varphi ,\nabla \varphi)$.\\
If $(f_k)_k$ is a sequence in $H_0^1(\Omega, \varphi ,\nabla \varphi)$ converging weakly to $0,$ then $(f_k)_k$ is also bounded in $L^2(\Omega,\varphi)$ and our assumption implies that we can find for any $N\in \mathbb N$ a smoothly bounded domain $\Omega_N\subset\subset\Omega$ such that
$$\Psi (z)=|\nabla \varphi (z)|^2+(1+\epsilon)\triangle \varphi (z)>N$$
on $\Omega\setminus\Omega_N$. Therefore we obtain
\begin{align*}
\int\limits_{\Omega}|f_k|^2e^{-\varphi}\,d\lambda  & \le
\int\limits_{\Omega_N}|f_k|^2e^{-\varphi }\,d\lambda  +
\int\limits_{\Omega\setminus\Omega_N} \frac{\Psi  |f_k|^2}{N}e^{-\varphi }\,d\lambda\\
&\le\|f_k\|^2_{L^2(\Omega_N)}+ \frac{C_\theta }{N} \|f_k\|^2_{1, \varphi, \nabla\varphi}.
\end{align*}
Now the classical Rellich -- Kondrachov Theorem asserts that the injection $H^1(\Omega_N) \hookrightarrow L^2(\Omega_N)$ is compact. Combined with our assumption, this shows that a subsequence of $(f_k)_k$ tends to $0$ in $L^2(\Omega, \varphi)$, which proves the Proposition.
\begin{flushright}
$\square$
\end{flushright}

\textbf{Remark. }Note that one does not need plurisubharmonicity of the weight function in the proof of the Proposition. If it is plurisubharmonic, one can of course drop $\theta$.
\vskip 0.5cm
Note also that interchanging the roles of $X_j$ and $X_j^*$ in the proof gives a criterion for compactness of the injection $H^1_0(\Omega,\varphi)\hookrightarrow L^2(\Omega,\varphi)$, since $\|f\|^2_{H_0^1(\Omega, \varphi)}=\| f\|^2_\varphi+\sum_{j=1}^n(\|X^*_jf\|_\varphi ^2+\|Y^*_jf\|_\varphi ^2)$. We formulate this in the next Proposition.

\begin{proposition} 
Suppose that the weight function satisfies
\begin{align*}
&\lim_{z\in\Omega, |z|\to \infty}(\theta |\nabla \varphi (z)|^2-\triangle \varphi (z))= +\infty \ \ \text{as well as}\\
&\lim_{z\in\Omega, z\to \partial\Omega}(\theta |\nabla \varphi (z)|^2-\triangle \varphi (z))= +\infty.
\end{align*}
for some $\theta\in (0,1)$.Then the embedding of $H_0^1(\Omega, \varphi)$ into $ L^2(\Omega, \varphi) $ is compact.
\end{proposition}
\vskip 0.5cm

\textbf{Remark.} The two above conditions are not sharp, which is not surprising since they do not take the geometry of the boundary into account. To see this, take $\varphi\equiv 0$, so both $H_0^1(\Omega, \varphi)$ and $H_0^1(\Omega, \varphi ,\nabla \varphi)$ coincide with the classical Sobolev space. But the injection $H_0^1(\Omega)\hookrightarrow L^2(\Omega)$ can be compact, if $\Omega$ is sufficiently thin at infinity. See \cite{adams}, Chapter 6 for various conditions.
\vskip 1 cm

\section{Compactness in the weighted problem.}~\\

The following Proposition is a well-known characterization of compactness in the $\dquer$-Neumann problem on bounded domains. In fact, it can be proven verbatim as for instance in \cite{str} in our context. 

\begin{proposition}
\label{compact}
Suppose that $\varphi$ is a plurisubharmonic weight function such that a bounded $\dquer$-Neumann operator $N_\varphi$ exists and let $\Vert.\Vert_X$ be a norm on $L^2_\varphi(\Omega)$ strictly weaker than $\Vert.\Vert_\varphi$. Then the following are equivalent:
\begin{enumerate}
\item The $\dquer$-Neumann operator on the level of $(0,1)$-forms $N_{1,\varphi}$ on $L^2_{(0,1)}(\Omega, \varphi)$  is compact.
\item The embedding of the space $dom(\dquer)\cap dom(\dquers)$ provided with the graph norm $f \mapsto(\Vert f\Vert_\varphi^2+\Vert \dquer f\Vert^2_\varphi+\Vert \dquers f\Vert^2_\varphi)^{\frac{1}{2}}$ into $L^2_{(0,1)}(\Omega,\varphi)$ is compact.
\item For each $\varepsilon> 0$ there exists a constant $C_\varepsilon>0$ such that
\begin{equation*}
\Vert f\Vert_\varphi\leq\varepsilon(\Vert \dquer f\Vert_\varphi^2+\Vert \dquers f\Vert^2_\varphi)^{\frac{1}{2}}+C_\varepsilon\Vert f\Vert_{X}
\end{equation*}
for all $f\in dom(\dquer)\cap dom(\dquers)$.
\end{enumerate}
\end{proposition}
To prove a first result on compactness of $N_\varphi $, we will make use of G\aa rding's inequality, which we now reformulate to suit in our context.

\begin{proposition}
\label{garding}
\textbf{(G\aa rding's inequality) } Let $\Omega$ be a smooth bounded domain. Then for any $f\in H^1(\Omega , \varphi ,\nabla \varphi)$ with compact support in $\Omega$,
\begin{equation*}
\Vert f \Vert ^2_{1,\varphi,\nabla\varphi}\leq C(\Omega, \varphi)\left(\Vert\dquer f\Vert^2_\varphi+\Vert\dquer^*_\varphi f\Vert^2_\varphi +\Vert f\Vert^2_\varphi\right) .
\end{equation*}
\end{proposition}

{\em Proof.}
For the proof we refer the reader to \cite{gh}, Propostion 4.3.
\begin{flushright}
$\square$
\end{flushright}

Following Catlin's idea for showing a sufficient condition for compactness of the $\dquer$-Neumann operator on bounded domains in \cite{Ca}, we prove the next Proposition. Indeed, we can use the same proof with only minor modifications, which arise from the fact that we are using a different norm.
\begin{proposition}
\label{cp 1}
Let $\Omega$ be a smooth pseudoconvex domain and let $\varphi$ be plurisubharmonic on $\Omega$. If the lowest eigenvalue $\lambda_\varphi (z)$ of the complex Hessian of $\varphi$ satisfies
\begin{equation}
\label{star}
\lim_{z\in\Omega, z\to\partial\Omega}\lambda_\varphi=\infty \quad\text{as well as}\quad\lim_{z\in\Omega, z\to\infty}\lambda_\varphi(z)=\infty,
\end{equation}
then $N_\varphi$ is compact.
\end{proposition}

{\em Proof.} By assumption and plurisubharmonicity of the weight we are in the setting of Proposition \ref{rellich}, thus it suffices to use Proposition \ref{compact} and show a compactness estimate.\\ 
Given $\epsilon>0$ we choose  $M\in\mathbb N$ with $1/M \le \epsilon/2$ and  a smooth bounded domain $\Omega_M\subset\subset\Omega$ such that $\lambda_\varphi(z)>M$ whenever $z\in\Omega\setminus\Omega_M$. Let $0\le\chi\le1$ be a smooth function with compact support in $\Omega$, which is identically one on $\Omega_M$. Hence we can estimate
\begin{align*}
M\Vert f\Vert_\varphi^2\leq& \sum_{j,k} \int\limits_{\Omega\setminus \Omega_M}\levim f_j\overline f_k e^{-\varphi}\,d\lambda +M\Vert\chi f\Vert_\varphi^2\\
\leq&Q_\varphi (f,f)+M\langle\chi f,f\rangle_\varphi\\
\leq&Q_\varphi (f,f)+M\Vert\chi f\Vert_{H_0^1(\Omega, \varphi ,\nabla \varphi)}\Vert f\Vert_{H_0^{-1}(\Omega, \varphi ,\nabla \varphi)}\\
\leq&Q_\varphi (f,f)+Ma\Vert\chi f\Vert_{1,\varphi,\nabla\varphi}^2+a^{-1}M\Vert f\Vert_{H_0^{-1}(\Omega, \varphi ,\nabla \varphi)}^2 ,
\end{align*}
where $a$ is to be chosen later. By assumption and Theorem \ref{existence}, a bounded $\dquer$-Neumann operator exists, which implies $\| f\|_\varphi^2\le C_\varphi (\| \dquer f\|^2_\varphi+\|\dquers f\|^2_\varphi)$ for some $C_\varphi>0$. Thus applying  G\aa rding's inequality \ref{garding} to the second  term, we find a constant $C_M$ only depending on $\Omega_M$, $\chi$ and $\varphi$ such that 
$$M\Vert f\Vert_\varphi^2\leq Q_\varphi (f,f)+MaC_MQ_\varphi (f,f)+a^{-1}M\Vert f\Vert_{H_0^{-1}(\Omega, \varphi ,\nabla \varphi)}^2.$$
Now choose $a$ such that $aC_M \le \epsilon/2,$ then
\begin{equation*}
\Vert f\Vert_\varphi^2\leq\epsilon Q_\varphi (f,f)+a^{-1}\Vert f\Vert_{H_0^{-1}(\Omega, \varphi ,\nabla \varphi)}^2
\end{equation*}
and this estimate implies compactness of $N_\varphi$ by Proposition \ref{compact}.
\begin{flushright}
$\square$
\end{flushright}

This condition on the weight function is of course rather restrictive, and it does not take the geometry of the boundary into account.  To weaken it, we first consider the following example.
\vskip 0.5cm

\textbf{Example. } Suppose that $\Omega\subset\mathbb C$ is the upper halfspace, given by $\Omega=\lbrace z:\mathfrak{Im}z>0\rbrace$. Let $\varphi_M= e^{-My}$, where $y=\mathfrak{Im}z$. Then, clearly, $0\le\varphi_M\le 1$ on $\Omega$ and $\varphi_M$ is subharmonic since $\triangle\varphi_M=M^2e^{-My}$. In particular $\triangle\varphi_M=M^2$ on  $\partial\Omega$. If we set
\begin{equation}
\label{phi bound}
\varphi =\sum_{j=1}^\infty \frac{1}{2^j}\varphi_{2^j},
\end{equation}
then $\varphi$ equals a bounded function on $\overline\Omega$ that is smooth and subharmonic in $\Omega$, such that $\triangle \varphi\to \infty$ as $z\to\partial \Omega$. This consideration shows that given a plurisubharmonic weight $\psi$ on $\Omega$ such that the lowest eigenvalue $\lambda_\psi$ of $M_\psi$ fulfills $\lambda_\psi\to\infty$ for $|z|\to\infty$, one can always construct by the substitution $\psi\mapsto\psi+\varphi$ a weight inducing an equivalent norm and satisfying the assumptions of Proposition \ref{cp 1}.\\
In a bit  more generality, suppose that $\Omega\subset\mathbb C^n$ admits a global defining function $r(z)$ which is strictly plurisubharmonic on $\overline\Omega$. Set again $\varphi_M(z)=e^{Mr(z)}$. As before, $0\le\varphi_M\le 1$ on $\Omega$ and $\varphi_M$ is strictly plurisubharmonic, since
$$\frac{\partial^2\varphi_M}{\partial z_j\partial\overline z_k}(z)=M\frac{\partial^2r}{\partial z_j\partial\overline z_k}(z)e^{Mr(z)}+M^2\frac{\partial r}{\partial z_j}(z)\frac{\partial r}{\partial\overline z_k}(z)e^{Mr(z)}.$$
Since we assumed strict plurisubharmonicity of $r$, the lowest eigenvalue of the Hessian $(\partial^2r/\partial z_j\partial\overline z_k)_{jk}$ is strictly positive for all $z\in\partial\Omega$. Although it could possibly tend to $0$ for $z\in\partial\Omega, |z|\to\infty$, a construction as in \eqref{phi bound} will nevertheless give us a bounded function such that the complex Hessian of $\varphi$ explodes at every boundary point, meaning that we proved the following Lemma.

\begin{lemma}
\label{explode boundary}
Let $\Omega\subset \mathbb C^n$ be an unbounded, smooth domain with strictly plurisubharmonic global defining function $r$. Then there is a bounded smooth plurisubharmonic function $\varphi$ on $\Omega$ such that all eigenvalues of the complex Hessian of $\varphi$ tend to infinity as $z\in\Omega$ tends to the boundary.
\end{lemma}

In particular we have the following Proposition.

\begin{proposition}
Let $\Omega\subset \mathbb C^n$ be an unbounded, smooth domain with a strictly plurisubharmonic global defining function $r$. Let $\lambda_\varphi$ be the lowest eigenvalue of the complex Hessian of the weight function. If
$$\lim_{z\in\Omega, |z|\to\infty}\lambda_\varphi=\infty,$$
then the $\dquer$-Neumann operator $N_\varphi$ is compact on $L^2_{(0,1)}(\Omega, \varphi)$.
\end{proposition}

{\em Proof. } By Lemma \ref{explode boundary} we can find an equivalent weight $\psi$ such that both 
$$\lim_{z\in\Omega, |z|\to\infty}\lambda_\psi=\infty\quad \text{and}\quad\lim_{z\in\Omega, z\to\partial\Omega}\lambda_\psi=\infty.$$
 By Proposition \ref{cp 1} $N_\psi$ is compact, thus also $N_\varphi$ since compactness is invariant under equivalent norms.
\begin{flushright}
$\square$
\end{flushright}

A similar construction also works under the weaker assumption that $\partial\Omega$ just satisfies Property (P), see Definition \ref{prop p}. Here we can not find a bounded function with properties as in Lemma \ref{explode boundary}, but nevertheless it is possible to construct for any given weight $\varphi$ with $\lambda_\varphi(z)\to\infty$ for $z\in \Omega, |z|\to\infty$ an equivalent one that fulfills condition \eqref{star}.\\
In order to proof Theorem \ref{cp 3}, we still need the following Proposition due to McNeal (see \cite{mcneal 2}, Proposition 2.1). We have to modify it a bit to suit our needs, so we also include the slightly changed proof.

\begin{proposition}
\label{mcneal}

\renewcommand{\theenumi}{\alph{enumi}}
Suppose that $\Omega$ is a smooth, unbounded and pseudoconvex domain in $\mathbb C^n$ and suppose that $V\subset\subset\mathbb C^n $ is open. If $V\cap\Omega\neq\emptyset$, then there exists a smooth, bounded and pseudoconvex domain $\tilde\Omega$ with the following properties:
\begin{enumerate}
\item $V\cap\Omega\subset \tilde\Omega$;
\item all points in $\partial\tilde\Omega\setminus\partial\Omega$ are strictly pseudoconvex.
\end{enumerate}
\renewcommand{\theenumi}{\arabic{enumi}}
\end{proposition}

{\em Proof.}
Let $K_1=\partial \Omega\cap V$. There is an interger $R$ such that $K_1\subset\subset K_2=\mathbb B_R\cap\partial\Omega$. Now $K_2$ is part of the boundary of a smooth bounded pseudoconvex domain $\tilde\Omega_2$, such that $V_1=\mathbb B_R\cap\Omega\subset \tilde\Omega_2$. To see this, intersect $\Omega$ with $\mathbb B_{R+1}$ to get a bounded pseudoconvex domain with continuous boundary and approximate it afterwards by something smooth. Thus we can use Theorem 1 from \cite{df} and find a smooth function $r$ defining $\tilde\Omega_2$, such that $-(-r)^\eta$ is strictly plurisubharmonic in a neighborhood of $K_2$ for some $0<\eta<1$. We can assume without loss of generality that  $-(-r)^\eta$ is strictly plurisubharmonic on $V_1$.\\
Choose $R_1$ such that $\mathbb B_{R_1}$ intersects $\partial \Omega$ transversally and such that $V\cap\Omega\subset \mathbb B_{R_1}\cap\Omega$. If we choose $V_1$ big enough, we can also assume $\mathbb B_{R_1}\cap\Omega\subset V_1$. So $dr(z)$ and $d|z|^2$ are linearly independent for $z\in\partial\mathbb B_{R_1}\cap\partial \Omega$, and by continuity this also holds on a neighborhood of the intersection. Thus we find $\rho_1,\rho_2$ and $\delta$, such that this is true on $(\mathbb B _{\rho_1}\setminus\mathbb B_{\rho_2})\cap\lbrace -\delta<r(z)<\delta\rbrace$.\\
From here on, we can follow verbatim McNeal's proof. Let $\chi_1(t)$ be a real-valued, smooth and increasing function, such that $\chi_1\equiv 0$ for $t\le R_1$ and $\chi_1^\prime(t)$ and $\chi_1^{\prime\prime}(t)$ strictly positive for $t>R_1$. Let $\chi_2(t)$ be smooth and increasing, such that $\chi_2\equiv -\delta^\eta$ for $t\le-\delta^\eta$ and $\chi_2(t)=t$ for $t\ge-\frac{1}{2}\delta^\eta$.\\
Now set 
$$\rho(z)=\chi_1(|z|^2)+\chi_2(-(-r(z))^\eta).$$
By the same calculation as in \cite{mcneal 2}, Propostion 2.1, one verifies that the domain defined by $\rho$ has the desired properties.
\begin{flushright}
$\square$
\end{flushright}

{\em Proof of Theorem \ref{cp 3}.}
First choose an arbitrary integer $M$. By assumption, one finds $R_0$ such that $\lambda _\varphi(z)>2^M$ for $|z|>R_0$. By Proposition \ref{mcneal} there is a smooth bounded pseudoconvex domain $\Omega_1$, such that $\Omega\cap\mathbb B_{R_0}\subset\Omega_1$ and $\partial\Omega_1\setminus\partial\Omega$ is strictly pseudoconvex. By assumption and strict pseudoconvexity of the rest of the boundary, $\Omega_1$ satisfies Property (P). So after choosing $M$ we can find $\varphi_1\in \mathcal C^\infty(\overline\Omega_1)$ with $0\le\varphi_1\le 1$ and  the lowest eigenvalue of the complex Hessian of $\varphi_1$ greater than $2^M$ on $\partial \Omega_1$. $\varphi_1$ is smooth on a closed set with smooth boundary, hence we can extend it smoothly to a bigger one. So extend $\varphi_1$ to a function $\psi_1\in \mathcal C^\infty(\overline\Omega)$, such that $0\le\psi_1\le 2$ and $\psi$ vanishes outside a ball with radius $R_1$. We can choose $R_1$ so big that the lowest eigenvalue of the complex Hessian of $\psi_1$ is bounded from below by $-2^{M-1}$ and $\lambda_\varphi (z)>2^{M+1}$  for $|z|>R_1$. Next we find a smooth bounded pseudoconvex domain $\Omega_2$ containing $\Omega\cap\mathbb B_{R_1}$, such that $\partial \Omega_2\setminus\partial\Omega$ is strictly pseudoconvex. $\Omega_2$ has Property (P). Hence there exists $\varphi_2$ with $0\le\varphi_2\le 1$ and lowest eigenvalue of the complex Hessian of $\varphi_2$ greater than $2^{M+1}$ on the boundary. Extend $\varphi_2$ to a function $0\le\psi_2\le2\in \mathcal C^\infty(\overline\Omega)$ with support in $\Omega\cap\mathbb B_{R_2}$ and Hessian bounded from below by $-2^M$ and $\lambda_\varphi(z)>2^{M+2}$ for $|z|>R_2$.\\
Inductively, we construct functions $\psi_j$ and by construction, $$\psi=\varphi+\sum\frac{1}{2^j}\psi_j$$ is a weight equivalent to $\varphi$ satisfying \eqref{star}. Therefore $N_\psi$ is compact by Proposition \ref{cp 1}, hence also $N_\varphi$.
\begin{flushright}
$\square$
\end{flushright}

{\bf Remark. }
Motivated by McNeal's generalization Property (\~ P) of Property (P) given in \cite{mc}, it would by interesting to know wether a version of Theorem \ref{cp 3} involving Property (\~ P) still holds true. Note that the proof given here heavily relies on the boundedness of the ``Property (P)''- functions.
\vskip 0.5 cm

\textbf{Example. } As was shown by Catlin in \cite{Ca}, all domains of finite type satisfy Property (P), so this provides a class of domains for which the Theorem can be applied. Consider for instance a domain of the form $\Omega_p=\lbrace (z^\prime ,z_n)\in \mathbb C^n\ |\ \mathfrak{Im}z_n>p(z^\prime)\rbrace$, where $p(z^\prime)$ is a plurisubharmonic function. In $\mathbb C^2$, such domains are always of finite type. 
\vskip 0.5 cm

In view of Theorem \ref{cp 3}, it is worth pointing out the following remark.
\vskip 0.5 cm

\textbf{Remark. }
Suppose that $\Omega_2\subset\Omega_1$. Let $\varphi_1$ be a weight function on $\Omega_1$ and let $\varphi_2$ be the restriction of $\varphi_1$ to $\Omega_2$. Then $L^2_{\varphi_2}(\Omega_2)$ is continuously embedded in $L^2_{\varphi_1}(\Omega_1)$. But note that $N_{\varphi_2}$ is not the restriction of $N_{\varphi_1}$ to $L^2_{\varphi_2}(\Omega_2)$. This is because $ker_{\Omega_2}(\dquer)$ is not embedded in $ker_{\Omega_1}(\dquer)$. In particular, compactness of $N_{\varphi_1}$ does not imply compactness of $N_{\varphi_2}$.
\vskip 0.5 cm

It remains to give the proof of Theorem \ref{ness}.
\vskip 0.5cm

{\em Proof of Theorem \ref{ness}.}
In the first step of the proof we show that if $N_\varphi$ is compact and if $(f_n)_{n=1}^\infty$ is a normed sequence weakly tending to zero, then $\langle \boxphi f_n,f_n\rangle_\varphi\to\infty$ for $n\to\infty$. (In fact, this property is equivalent to compactness).\\
So let $N_\varphi$ be compact. By the Spectral Theorem for compact self-adjoint operators, there is an orthonormal basis of $L^2_{(0,1)}(\Omega,\varphi)$ consisting of eigenvectors of $N_\varphi$, call it $\lbrace v_j\rbrace_{j\in\mathbb N}$. We have $N_\varphi v_j=\lambda_j v_j$, where $\lambda_j\to0$ for $j\to\infty$ and we assume the $\lambda_j$ to be ordered decreasingly. Moreover we have $v_j\in dom(\boxphi)$ and $\boxphi v_j=1/\lambda_j v_j$. Now if $(f_n)_n$ is a normed sequence weakly converging to zero, then $f_n=\sum_{j=1}^\infty a_{nj}v_j$, where for all $n$ it holds $\sum_{j=1}^\infty |a_{nj}|^2=1$ and for all $j$ it holds that $|a_{nj}|\to0$ as $n\to\infty$, since weak convergence is equivalent to coordinatewise convergence. Hence, for any given $M\in\mathbb N$ and $\varepsilon>0$ we find $J$ such that $1/\lambda_j>M$ for all $j>J$ and after that $N$ such that $\sum_{j=1}^J |a_{nj}|^2<\varepsilon$ for any $n>N$. Thus for any $n>N$
\begin{align*}
\langle \boxphi f_n,f_n\rangle_\varphi=&\sum_{j=1}^J\frac{1}{\lambda_j}|a_{nj}|^2+\sum_{j=J+1}^\infty\frac{1}{\lambda_j}|a_{nj}|^2\\
\ge&M(1-\varepsilon)\\
>&M/2
\end{align*}
for $\varepsilon$ sufficiently small, which proves the first statement (Note that in this computation, one can not directly commute $\square _\varphi$ with the infinite sum, since it is not bounded. Nevertheless the identity holds true, as one can see be substituting $f_n=N_\varphi u_n$ and using the uniqueness of the expression in an orthonormal basis).\\
To finish the proof, let $(\mathbb B (z_l,r))_l$ be a sequence of disjoint balls in $\Omega$. Without loss of generality we can assume $z_0=0$. Now let $\chi\in\Lambda _0^{(0,1)}(\mathbb B (0,r))$ be a real-valued form with $\| \chi\|_{L^2_{(0,1)}}=1$ and set $\chi^{(l)}=\chi (z-z_l)$ as well as $f^{(l)}=\chi^{(l)}e^{\varphi/2}$. Then $(f^{(l)})_l$ is a normed sequence in $L^2_{(0,1)}(\Omega,\varphi)$,  and it tends weakly to zero since the support of $f^{(l)}$ moves out to infinity. Thus by the first part of the proof combined with the Kohn -- Morrey formula \ref{kohn-morrey} we have that
$$\sum_{j,k=1}^{n}\ \int \limits_{\Omega}\levim f^{(l)}_j\overline f^{(l)}_ke^{-\varphi}d\lambda+\sum_{j,k=1}^{n}\ \int \limits_{\Omega}\left| \frac{\partial f^{(l)}_j}{\partial \overline z_k}\right| ^2e^{-\varphi}d\lambda\to\infty$$
as $l\to\infty$. Plugging in the definition of $f^{(l)}$, we get
$$\sum_{j,k=1}^{n}\ \int \limits_{\mathbb B (z_l,r)}\levim \chi^{(l)}_j \chi^{(l)}_kd\lambda+\sum_{j,k=1}^{n}\ \int \limits_{\mathbb B (z_l,r)}\left| \chi^{(l)}_j\frac{\partial \varphi}{\partial \overline z_k}+\frac{\partial\chi^{(l)}_j}{\partial\overline z_k}\right| ^2d\lambda\to\infty.$$
Now $\chi^{(l)}$ and its first order derivatives are uniformly bounded in $l$, so estimating the complex Hessian of $\varphi$ by its trace and using the triangle inequality it follows
$$\lim_{l\to\infty}\ \sum_{j=1}^n\int_{\mathbb B (z_l,r)}|\chi^{(l)}|^2 \left(\triangle\varphi+|\nabla\varphi|^2\right)d\lambda=\infty.$$
After increasing $r$, this combined with Lemma \ref{vergl} implies the claim.
\begin{flushright}
$\square$
\end{flushright}

{\bf Remark. }Note that by Lemma \ref{vergl}, $\lim_{l\to\infty}\int_{\mathbb B (z_l,r)}\triangle\varphi\ d\lambda=\infty$ for any sequence $\mathbb B (z_l,r)$  of disjoint balls  relatively compact in $\Omega$ holds if and only if 
$$\lim_{l\to\infty}\int_{\mathbb B (z_l,r)}|\nabla\varphi|^2\ d\lambda=\infty$$
for any such sequence.
\vskip 1cm

\section{Some applications to the unweighted problem.}~\\
\label{unweighted}

In this closing section, we will think of the unweighted problem as the special case $\varphi\equiv 0$. Let us start with an illustrative example of what one can expect.
\vskip 0.5 cm

\textbf{Example. } Suppose that $\Omega$ is not quasibounded, so there exists a sequence of disjoint balls $\mathbb B(z_l,r)$ with fixed radius $r$ contained in $\Omega$. Consider a $(0, 1)$-form $v\in\Lambda^{(0,1)}_0(\mathbb B (z_0, r))$ such that $\| v\|_{L^2_{(0,1)}} = 1$. Without loss of generality we can assume $z_0 = 0$ and clearly, $v\in dom( \square )$. Now take transverses $v_n(z) = v(z-z_n)$ of $v$. By definition, they have disjoint support and $\|v_n\|_{L^2_{(0,1)}} = 1$. $(\square v_n)_n$ is a bounded sequence, but the functions $N\square v_n= v_n$ are pairwise orthogonal hence they can not contain a convergent subsequence. Thus N is not compact. Note also that $\varphi\equiv 0$ does not satisfy the necessary condition of Theorem \ref{ness}.\\
Suppose that there is a sequence $\mathbb B(z_l, r_l)$ of disjoint balls contained in $\Omega$, such that $r_l\to\infty$, then by a similar argument it  follows that the unweighted $\dquer$-Neumann operator on $\Omega$ is not bounded.
\vskip 0.5cm

\textbf{Remark. }The unit ball $\mathbb B\subset\mathbb C^n$ is strictly pseudoconvex and thus the $\dquer$-Neumann operator on $\mathbb B$ compact. Nevertheless, the ball is biholomorphic to the Siegel upper half space $\mathbb U=\lbrace (z^\prime,z_n)\in\mathbb C^n: \mathfrak{Im}z_n>|z^\prime|^2\rbrace$ via the Cayley transform, and the previous example shows that there is no bounded $\dquer$-Neumann operator on $\mathbb U$. This shows in particular, that existence and compactness in the $\dquer$-Neumann problem are not invariant under biholomorphisms. 
\vskip 0.5cm

The example also motivates the following Definition, which is again taken from \cite{adams}.

\begin{definition}
A domain $\Omega$ is called quasicylindrical if and only if $$\limsup_{z\in\Omega, \ |z|\to\infty} dist(z,\partial\Omega)\le C$$
for some $C>0$. Equivalently, $\Omega$ is quasicylindrical if and only if there is no sequence of pairwise disjoint balls with radii going to infinity contained in $\Omega$.
\end{definition}

Summing up, we can state the following Lemma.

\begin{lemma}
Suppose that $\Omega$ is an unbounded domain. If the there is a bounded $\dquer$-Neumann operator on $\Omega$, then $\Omega $ is quasicylindrical. If the $\dquer$-Neumann operator is compact, then $\Omega$ is quasibounded.
\end{lemma}

On the other hand, we can combine the Kohn -- Morrey formula \ref{kohn-morrey} with the fact that existence of the $\dquer$-Neumann operator is invariant under equivalent weights, to get a sufficient condition for existence.

\begin{lemma}
\label{bounded}
Suppose that $\Omega$ is pseudoconvex and that there exists a bounded plurisubharmonic function $\varphi$ on $\Omega$,  such that 
$$\liminf_ {z\in\Omega, |z|\to\infty}\lambda_\varphi\ge\varepsilon>0.$$
Then there exists a bounded $\dquer$-Neumann operator on $L^2(\Omega)$.
\end{lemma}

{\em Proof.}
If $\varphi$ is a function with the assumed properties, consider the $\dquer$-Neumann problem in the weighted space $L^2_{(0,1)}(\Omega,\varphi)$. By Theorem \ref{existence}, there is a bounded $\dquer$-Neumann operator on $L^2_{(0,1)}(\Omega,\varphi)$, hence also on $L^2_{(0,1)}(\Omega)$, since a bounded weight is equivalent to the one identically zero.
\begin{flushright}
$\square$
\end{flushright}

\textbf{Example. } Let $\Omega$ be of the form $\Omega=D\times \lbrace-1<\mathfrak {Im}z_n<1\rbrace$, where $D$ is a bounded domain in $\mathbb C^{n-1}$. Then $\varphi=\|z^\prime\|^2+y_n^2$ is a bounded plurisubharmonic function on $\Omega$, such that $\lambda_\varphi\ge\frac{1}{2}$.\\
Conversely, suppose that $\Omega$ contains a complex line. Then there can be no such function, since there is no bounded plurisubharmonic function $\varphi$ on $\mathbb C$ such that $\triangle\varphi>\varepsilon$ uniformly.
\vskip 0.5cm

\textbf{Remark. }The same argument can be used to show existence of a $\dquer$-Neumann operator on non-pseudoconvex domains, if one assumes the Levi-form of the defining function of $\Omega$ to be semibounded from below. Suppose that 
$$\sum_{j,k=1}^n\frac{\partial^2r}{\partial z_j\partial\overline z_k}(z)t_j\overline t_k\ge -C\|t\|^2\quad \text{for all }t\in\mathbb C^n,\ z\in\Omega$$
and assume that there is a plurisubharmonic function $\varphi$ on $\Omega$ as in Lemma \ref{bounded}. Then $\psi=(C+1)\varphi/\varepsilon$ is a bounded function and the Kohn -- Morrey formula \ref{kohn-morrey} shows that $Q_\psi(f,f)\ge\| f\|_\psi^2$, thus a bounded $\dquer$-Neumann operator exists.\\
In particular if $\Omega$ is bounded and of class $\mathcal C^2$, there is a $\mathcal C^2$-defining function $r$ and since $\partial \Omega$ is compact, the Levi-form of $r$ is always bounded from below. $\varphi(z)=\| z\|^2$ is a bounded plurisubharmonic function on $\Omega,$ with $\lambda_\varphi=1$, thus a bounded $\dquer$-Neumann operator exists on each bounded domain with $\mathcal C^2$-boundary.
\vskip 0.5cm

\begin{lemma}
\label{cp 4}
Let $\Omega$ be pseudoconvex and suppose that it satisfies Property (P). Suppose furthermore that there is a bounded plurisubharmonic function $\varphi$ on $\Omega$, such that $\lambda_\varphi\to\infty$ for $|z|\to\infty$. Then $N$ is compact on $L^2(\Omega)$.
\end{lemma}

{\em Proof.} Take $\varphi$ as weight function. Theorem \ref{cp 3} assures that $N_\varphi$ is compact on $L^2_{(0,1)}(\Omega,\varphi)$, thus also the unweighted $\dquer$-Neumann operator on $\Omega$ is compact, since compactness is invariant under equivalent weights.
\begin{flushright}
$\square$
\end{flushright}

\textbf{Example. } Suppose that $\Omega$ is given by $\Omega=\lbrace x+iy\in\mathbb C\ |\ x^2y^2<1\rbrace$. Then by definition, $\varphi=x^2y^2$ is a bounded function on $\Omega$ and furthermore we have $\triangle\varphi=x^2+y^2$. Thus $N$ is compact on $L^2(\Omega)$ by Lemma \ref{cp 4}.
\vskip 0.5cm

\textbf{Acknowledgement. } The author would like to thank Anne-Katrin Herbig for many useful and stimulating discussions on this topic.

\end{document}